\documentclass[12pt]{article}
\usepackage{bbm}
\usepackage{latexsym}
\usepackage{mathrsfs}
\usepackage{graphicx}
\usepackage{subfigure}
\usepackage{amssymb}
\usepackage{amsmath}
\usepackage{amsthm}
\usepackage{color}
\usepackage{cite}
\usepackage{float}
\usepackage{indentfirst}
\usepackage{anysize}\marginsize{45mm}{45mm}{40mm}{50mm}
\usepackage{caption}
\usepackage{float}
\usepackage{multirow}

\usepackage{enumerate}

\usepackage{latexsym, bm}
\newtheoremstyle{lemma}{\topsep}{\topsep}%
     {}
     {}
     {\bfseries}
     {}
     {0.1em}
     {\thmname{#1}\thmnumber{ #2}\thmnote{ #3}}
\theoremstyle{lemma}  

\newtheorem{theorem}{Theorem}     
\newtheorem{lemma}[theorem]{Lemma}
\newtheorem{corollary}[theorem]{Corollary}
\newtheorem{conjecture}[theorem]{Conjecture}

\numberwithin{equation}{section}

\title{ On the conjecture of bijection between perfect matching and sub-hypercube in folded hypercubes \thanks{This research was partially supported by the National Natural Science Foundation of China (Nos. 11801061 and 11761056), the Chunhui Project of Ministry of Education (No. Z2017047) and the Fundamental Research Funds for the Central Universities (No. ZYGX2018J083)}}

\author{ Huazhong L\"{u}$^{1}$\thanks{Corresponding author.} and Tingzeng Wu$^{2}$\\
{\small $^{1}$School of Mathematical Sciences, University of Electronic Science and Technology of China,} \\
{\small Chengdu, Sichuan 610054, P.R. China}\\
{\small E-mail: lvhz08@lzu.edu.cn}\\
{\small $^{2}$School of Mathematics and Statistics, Qinghai Nationalities University, }\\
{\small Xining, Qinghai 810007, P.R. China} \\
{\small E-mail: mathtzwu@163.com}\\}
\date{}
\begin{document}

\maketitle
\begin{abstract}

Dong and Wang in [Theor. Comput. Sci. 771 (2019) 93--98] conjectured that the resulting graph of the $n$-dimensional folded hypercube $FQ_n$ by deleting any perfect matching is isomorphic to the hypercube $Q_n$. In this paper, we show that the conjecture holds when $n=2,3$, and it is not true for $n\geq4$.

\vskip 0.1 in


\end{abstract}

\section{Introduction}

Let $G=(V(G),E(G))$ be a graph, where $V(G)$ is the vertex-set of $G$ and $E(G)$ is the edge-set of $G$. A {\em matching} of $G$ is a set of pairwise nonadjacent edges. A {\em perfect matching} of $G$ is a matching with size $|V(G)|/2$. A {\em $k$-factor} of $G$ is a $k$-regular spanning subgraph of $G$. Clearly, a perfect matching, together with end vertices of its edges, forms a 1-factor of $G$. $G$ is {\em $k$-factorable} if it admits a decomposition into $k$-factors. The distance between two vertices $u$ and $v$ is the number of edges in a shortest path joining $u$ and $v$ in $G$, denoted by $d(u, v)$. For any two edges $uv$ and $xy$, the distance of $uv$ and $xy$, denoted by $d(uv, xy)$, is $\min\{d(u, x), d(u, y), d(v, x), d(v, y)\}$. For other standard graph notations not defined here please refer to \cite{Bondy}.

\vskip 0.0 in

\begin{figure}
\centering
\includegraphics[height=40mm]{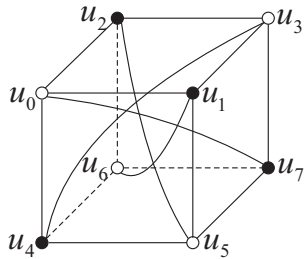}
\caption{The 3-dimensional folded hypercube $FQ_3$.} \label{FQ3}
\end{figure}

The well-known $n$-dimensional hypercube is a graph $Q_n$ with $2^n$ vertices and $n2^{n-1}$ edges. Each vertex is labelled by an $n$-bit binary string. Two vertices are adjacent if their binary string differ in exactly one bit position. The folded hypercube, denoted by $FQ_n$, is first introduced by El-Amawy and Latifi\cite{El-Amawy} as a variant of the hypercube. $FQ_n$ is obtained from the hypercube $Q_n$ by adding $2^{n-1}$ independent edges, called {\em complementary edges}, each of which is between $x_1x_2\cdots x_n$ and $\overline{x}_1\overline{x}_2\cdots \overline{x}_n$, where $\overline{x}_i=1-x_i$, $i=1,\cdots,n$. For convenience, the set of complementary edges of $FQ_n$ are denoted by $E_c$ and the set of $i$-dimensional edges in $Q_n$ are denoted by $E^i$ for each $1\leq i\leq n$, where an edge $uv$ is $i$-dimensional in $Q_n$ if $u$ and $v$ differ only in the $i$-th position. We illustrate $FQ_{3}$ in Fig. \ref{FQ3}.

Some attractive properties of the folded hypercube are widely studied in the literature, such as, pancyclicity \cite{Xu}, conditional connectivity \cite{Zhao}, stochastic edge-fault-tolerant routing algorithm \cite{Thuan}, conditional diagnosability \cite{Liu} and conditional cycle embedding\cite{Kuo}. Recently, Dong and Wang \cite{Dong} conjectured the following:

\begin{conjecture}{\bf.}
An subset $E^m$ of $2^{n-1}$ edges of $FQ_n$ is a perfect matching if and only if $FQ_n-E^m$ is isomorphic to $Q_n$.
\end{conjecture}

\vskip 0.05 in

We solve this conjecture in Section 2. Conclusions are given in Section 3.

\section{Main results}
%

%
%

The affirmative answer to Dong's conjecture for $n=2,3$ are shown as follows.

\begin{theorem}\label{iso}{\bf.} For any perfect matching $M$ of $FQ_n$, $n=2,3$, $FQ_n-M$ is isomorphic to $Q_n$.
\end{theorem}
\noindent{\bf Proof.} Clearly, $FQ_2$ is the complete graph $K_4$ and $Q_2$ is a 4-cycle, so the statement holds when $n=2$. Let $M$ be a perfect matching of $FQ_3$. If $M=E_c$, the lemma is obviously true. Therefore, we assume that $M\neq E_c$. For convenience, we label each vertex of $FQ_3$ by $u_i$, $i\in\{0,\cdots,7\}$, respectively (see Fig. \ref{FQ3}). We distinguish the following cases.

\noindent{\bf Case 1.} $M$ contains no complementary edges. By symmetry of $FQ_3$, we may assume that $u_0u_2\in M$. Then the subgraph induced by the hypercube edges of $FQ_3-\{u_0,u_2\}$, say $H$, is a $2\times 1$-grid with six vertices. Clearly, $H$ has three perfect matchings $M_1=\{u_1u_5,u_3u_7,u_4u_6\}$,  $M_2=\{u_1u_3,u_4u_5,u_6u_7\}$ and $M_3=\{u_1u_3,u_4u_6,u_5u_7\}$. Thus, $M=\{u_0u_2\}\cup M_j$, $j=1,2,3$. By direct checking, $FQ_3-M$ is isomorphic to $Q_3$.

\noindent{\bf Case 2.} $M$ contains exactly one complementary edge. Suppose w.l.o.g. that $u_0u_7\in M$. Then the subgraph induced by the hypercube edges of $FQ_3-\{u_0,u_7\}$, say $C$, is a 6-cycle. Clearly, $C$ has two perfect matchings $M_1=\{u_1u_5,u_2u_3,u_4u_6\}$ and $M_2=\{u_1u_3,u_2u_6,u_4u_5\}$. Thus, $M=\{u_0u_7\}\cup M_1$ or $\{u_0u_7\}\cup M_2$. By direct checking, $FQ_3-M$ is isomorphic to $Q_3$.

\noindent{\bf Case 3.} $M$ contains exactly two complementary edges. By symmetry of $FQ_3$, we may assume that $A\subset M$ or $B\subset M$, where $A=\{u_0u_7,u_2u_5\}$ and $B=\{u_0u_7,u_3u_4\}$.
Then $A$ (resp. $B$) can be uniquely extended to a perfect matching of $FQ_3$ by adding two hypercube edges. Thus, $M=\{u_0u_7,u_2u_5,u_1u_3,u_4u_6\}$ or $\{u_0u_7,u_3u_4,u_1u_5,u_2u_6\}$. By direct checking, $FQ_3-M$ is isomorphic to $Q_3$.

\noindent{\bf Case 4.} $M$ contains exactly three complementary edges. In this condition, exactly six vertices of $FQ_3$ are saturated by the complementary edges of $M$ and the remaining two vertices are diagonal, which can not be saturated by any hypercube edge. So there exist no perfect matchings containing exactly three complementary edges. This completes the proof.\qed
\vskip 0.1 in

The following lemma is useful.

\begin{lemma}\cite{Zhu}\label{common-neighbors}{\bf.} Any two vertices in $V(FQ_n)$ exactly have two common neighbors for $n\geq 4$ if they have.
\end{lemma}

For $n\geq4$, in fact, we prove the following theorem which characterizes the relationship between a perfect matching and the sub-hypercube of $FQ_n$.

\begin{theorem}{\bf.}\label{non-iso} Let $n\geq4$ be an integer and let $M$ be a perfect matching of $FQ_{n}$. Then $FQ_{n}-M$ is isomorphic to $Q_n$ if and only if $M=E_c$ or $E^i$ for any $i\in\{1,2,\cdots,n\}$.
\end{theorem}

\noindent{\bf Proof.} \textit{Sufficiency.} By the definition of $FQ_n$, if $M=E_c$, the statement is obviously true. Therefore, let $M=E^i$ for some $i\in\{1,2,\cdots,n\}$. We shall show that $FQ_n-E^i$ is isomorphic to $Q_n$. One may consider the graph $FQ_n-E^i\cup E_c$ since $FQ_n-E^i\cup E_c$ is two disjoint copies of $Q_{n-1}$. For convenience, let $G=FQ_n-E^i$. The vertices of $G$ are still labelled by $n$-tuple binary strings. We define a bijection $\varphi: V(G)\rightarrow V(Q_n)$ as follows: (1) $\varphi(u)=u$ if the $i$-th bit of $u$ is 0; (2) $\varphi(u)=\overline{u}_1\cdots\overline{u}_{i-1}u_i\overline{u}_{i+1}\cdots \overline{u}_n$ if the $i$-th bit of $u$ is 1, where $u=u_1\cdots u_{i-1}u_iu_{i+1}\cdots u_n$. Let $uv\in E(G)$ be an arbitrary edge. We shall verify that $\varphi$ is an isomorphism.

\noindent{\bf Case 1.} $uv$ is a $j$-dimensional edge of $G$, $j\in\{1,\cdots,n\}\setminus\{i\}$. Then $u$ and $v$ differ only in the $j$-th position. We may assume that $u=u_1\cdots u_i\cdots u_j\cdots u_n$ and $v=u_1\cdots u_i\cdots \overline{u}_j\cdots u_n$. If $u_i=0$, then $\varphi(u)=u$ and $\varphi(v)=v$, yielding that $\varphi(u)\varphi(v)\in E(Q_n)$. If $u_i=1$, then $\varphi(u)=\overline{u}_1\cdots\overline{u}_{i-1} u_i\overline{u}_{i+1}\cdots \overline{u}_j\cdots \overline{u}_n$ and $\varphi(v)=\overline{u}_1\cdots\overline{u}_{i-1} u_i\overline{u}_{i+1}\cdots u_j\cdots \overline{u}_n$. Again, $\varphi(u)\varphi(v)\in E(Q_n)$.

\noindent{\bf Case 2.} $uv\in E_c$. For convenience, let $u=u_1\cdots u_i\cdots u_n$ and $v=\overline{u}_1\cdots \overline{u}_i$ $\cdots \overline{u}_n$. We may assume that $u_i=0$. Thus, $\varphi(u)=u$ and $\varphi(v)=u_1\cdots \overline{u}_i\cdots u_n$. Therefore, $\varphi(u)\varphi(v)\in E(Q_n)$. By above, it follows that $FQ_{n}-M$ is isomorphic to $Q_n$.

\textit{Necessity.} Suppose on the contrary that $M\neq E_c$ and $M\neq E^i$ for each $i\in\{1,2,\cdots,n\}$. We consider the following two cases.

\noindent{\bf Case 1.} $M\cap E_c\neq\emptyset$. We claim that there exists a vertex $u$ such that the complementary edge $uv\in E_c$ and one of its neighbors, say $v_1$, is saturated by a hypercube edge $v_1u_1$ in $M$. Suppose not. If all the neighbors of any vertex $u$ in $FQ_n$ are saturated by complementary edges, then $M=E_c$. So the claim holds. Thus, there exists a 4-cycle $uv_1u_1v_2u$ in $FQ_n$, where $v_1$ and $v_2$ are two neighbors of $u$. Obviously, $uv_1,uv_2\not\in M$. Note that $u_1v_1\in M$, then $u_1v_2\not\in M$. This implies that $u$ and $u_1$ have exactly one common neighbor in $FQ_n-M$, contradicting the well-known fact that every two vertices in $Q_n$ have zero or exactly two common neighbors.

\noindent{\bf Case 2.} $M\cap E_c=\emptyset$. Our objective is to show that there exists a 4-cycle $C$ of $FQ_n$ containing exactly one edge of $M$. Accordingly, two diagonal vertices of $C$ have exactly one common neighbor, which contradicts the fact that every two vertices in $Q_n$ have zero or exactly two common neighbors. Note that $M\cap E_c=\emptyset$ and $M\neq E^i$ for each $i\in\{1,\cdots,n\}$, then there exists two edges $e,f\in M$ with $e\in E^i$ and $f\in E^j$ such that $d_{FQ_n}(e,f)=1$, where $1\leq i<j\leq n$. For clarity, let $C=uxvyu$ and $e=ux$. Suppose w.l.o.g. that all edges of $C$ are hypercube edges. So there exists an edge $xw$ connecting $e$ and $f$, where $w$ is an end vertex of $f$.

By Lemma \ref{common-neighbors}, $u$ and $v$ have exactly two common neighbors $x$ and $y$ in $FQ_n$, and vice versa. Note that $e\not\in E(FQ_n-M)$, $u$ and $v$ have at most one common neighbor in $FQ_n-M$. If $u$ and $v$ have exactly one common neighbor, say $y$, then we are done. So we assume that $u$ and $v$ have no common neighbors in $FQ_n-M$, namely $vy\in M$.

\begin{figure}
\centering
\includegraphics[height=40mm]{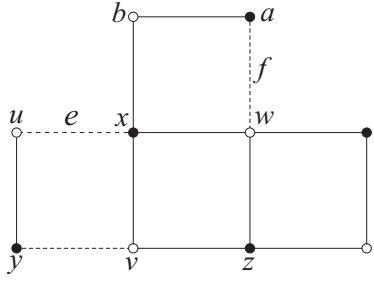}
\caption{Illustration for Theorem \ref{non-iso}.} \label{FQn-M}
\end{figure}

If $xv,uy\in E^j$, then there exists a 4-cycle $C'=xwzvx$ such that $f=wz$. Clearly, $xw,zv,vx\not\in M$ and $wz\in M$, then we have a 4-cycle that contains exactly one edge in $M$, yielding that $x$ and $z$ have exactly one common neighbor. So we assume that $xv,uy\not\in E^j$ and $f=wa\in M$. Accordingly, $xw\in E^k$, where $k\neq i,j$. Thus, there exists a cycle $C''=xwabx$ in $FQ_n$. Recall that $e=ux$ and $e\in M$, thus, $bx\not\in M$. Similarly, $ab,xw\not\in M$. So $a$ and $x$ have exactly one common neighbor in $G$, a contradiction (see Fig. \ref{FQn-M}). Hence, the theorem holds. \qed

By the above theorem, we have the following corollary, which disproves Dong's conjecture for $n\geq4$.

\begin{corollary}{\bf.} There exists a perfect matching $M$ of $FQ_n$ with $n\geq4$ such that $FQ_n-M$ is not isomorphic to $Q_n$.
\end{corollary}

\noindent{\bf Proof.} Obviously, there exists a perfect matching $M$ of $Q_n$ such that $M\neq E^i$ for each $1\leq i\leq n$. Note that $Q_n$ is a spanning subgraph of $FQ_n$, then $M\neq E_c$. By Theorem \ref{non-iso}, the statement follows immediately.\qed

\section{Conclusions}

In this paper, we characterize the relationship between the resulting graph of $FQ_n$ by deleting a perfect matching and the sub-hypercube. It is interesting to study the similar property in hypercube variants which include the hypercube as their spanning subgraphs.
%

\end{document}